\documentclass[a4study,twoside,11pt]{article}
\usepackage{amssymb}
\usepackage{amsmath}
\usepackage{graphicx}
\usepackage{amsthm}
\usepackage{cite}
\usepackage{array}

\usepackage{amsmath,latexsym,amssymb,amsfonts,amsbsy,amsthm}
\usepackage{graphicx}
\usepackage{subfigure}
\usepackage[a4paper]{geometry}
\usepackage{epstopdf}
\usepackage{diagbox}
\usepackage{enumerate}
\usepackage{color}

\usepackage{comment}
\usepackage{indentfirst}

\newfont{\bb}{msbm10}

\newtheorem{theorem}{Theorem}[section]
\newtheorem{definition}{Definition}[section]
\newtheorem{lemma}{Lemma}[section]

\newtheorem{remark}{Remark}[section]

\usepackage{verbatim}
\usepackage{algorithm}
\usepackage{algorithmicx}
\usepackage{algpseudocode}
\usepackage{amsmath}
\usepackage{graphics}
\usepackage{epsfig}
\usepackage{arydshln}

\usepackage{mathrsfs}
\usepackage{graphicx}
\usepackage{subfigure}
\usepackage{caption}
\usepackage{multirow}
\usepackage{threeparttable}
\numberwithin{equation}{section}

\usepackage{float}
\usepackage{footnote}
\makesavenoteenv{table}
\usepackage{booktabs}

\usepackage{authblk}

\usepackage{bm}
\usepackage{hyperref} 

\newenvironment{Proof}{{\noindent\it Proof.}}{\hfill $\square$ \par}

\title{On averaging block Kaczmarz methods for solving nonlinear systems of equations}
\author{
	A-Qin Xiao\\
	School of Mathematical Sciences, Tongji University,\\
	Shanghai, 200092, PR China.\\
	Email:xiaoaqin@tongji.edu.cn\\
        and\\
	Jun-Feng Yin\thanks{Corresponding author. This work was supported by National Natural Science Foundation of China (Grant No. 11971354)}\\
	School of Mathematical Sciences, Tongji University,\\
	Shanghai, 200092, PR China.\\
	Email:yinjf@tongji.edu.cn\\
}
\date{ }

\begin{document}
\maketitle
\begin{abstract}
A class of averaging block nonlinear Kaczmarz methods is developed for the solution of the nonlinear system of equations. 
The convergence theory of the proposed method is established under suitable assumptions and the upper bounds of the convergence rate for the proposed method with both constant stepsize and adaptive stepsize are derived.
Numerical experiments are presented to verify the efficiency of the proposed method, which outperforms the existing nonlinear Kaczmarz methods in terms of the number of iteration steps and computational costs.
 \\

\end{abstract}

\section{Introduction}\label{sec1}
Consider the solution of the nonlinear system of equations
\begin{equation}\label{eqnesys}
F(x)=0, 
\end{equation}
where $x\in \mathbb{R}^n$ is an unknown vector and $F(x)=(F_1(x),F_2(x),\cdots,F_m(x))^T$ is a nonlinear vector-valued function, which often arises in many scientific and engineering computations, for instance, optimization problems \cite{1975M,2011Y,2018H}, phase retrieval \cite{2015CLS,2020FS} and deep learning \cite{2016K,2022LZB}. 

The nonlinear Kaczmarz method \cite{2022YLG,2022WLB} is a simple and effective iterative method for solving the nonlinear system \eqref{eqnesys}. Let $F'(x)$ be the Jacobian matrix of $F(x)$ and $\nabla F_i(x)$ is the $i$-th row of matrix $F'(x)$. The nonlinear Kaczmarz method iterates by
\begin{equation}\label{eqiternkm}
    x_{k+1}= x_{k}-\frac{F_{i_k}(x_k)}{ \left\|\nabla F_{i_k}(x_k)\right\| ^2_2} \nabla F_{i_k}(x_k),
\end{equation}
where the index $i_k$ is cyclically or randomly selected from $[m]=\{1,\cdots,m\}$.  
Moreover, by choosing the index corresponding the maximum residual of partial or complete nonlinear equations, the nonlinear sampling Kaczmarz-Motzkin method \cite{2023ZBLW} and the maximum residual nonlinear Kaczmarz method \cite{2023ZWZ} were presented, respectively, and both of them accelerated the convergence rate of the existing nonlinear Kaczmarz methods. For more related studies on the nonlinear Kaczmarz methods, we refer the readers to \cite{2004MNT,2007DB,2007HLS}.

The block nonlinear Kaczmarz method which utilizes multiple nonlinear equations simultaneously at each iteration, can be implemented more efficiently and it updates by
\begin{equation}\label{eqitermrbnk}
    x_{k+1}= x_{k}- (F'_{\mathcal{I}_k}(x_k) )^{\dagger}F_{\mathcal{I}_k}(x_k),
\end{equation}
where $(F'_{\mathcal{I}_k}(x_k))^{\dagger}$ represents the Moore-Penrose pesudoinverse of the submatrix of $F'(x_k)$ and the block index set $\mathcal{I}_k\subseteq [m]$. 
In particular, the classical Newton-Raphson method \cite{1995Y,2022P} can be regarded as a special case of the block nonlinear Kaczmarz method when $\mathcal{I}_k= [m]$.
Moreover, a class of block sampling nonlinear Kaczmarz methods \cite{2022ZLT} were presented by randomly choosing a sample of $[m]$ to construct the index subset.
A maximum residual block nonlinear Kaczmarz method \cite{2023ZWZ} was proposed with an approximate greedy strategy to grasp almost the larger entries of the residual vector at each iteration. 

However, the above block nonlinear Kaczmarz methods require calculating the Moore-Penrose pseudoinverse of the selected submatrix from the Jacobian matrix, which usually costs expensive.
In this paper, to improve the efficiency of the block nonlinear Kaczmarz methods, a class of averaging block nonlinear Kaczmarz methods is developed for the solution of nonlinear systems of equations.
Theoretical analysis gives the upper bound on the convergence rate of the proposed method.
Numerical experiments show that the averaging block nonlinear Kaczmarz method is effective and faster than existing nonlinear Kaczmarz methods in terms of the number of iteration steps and CPU time.

The rest of this paper is organized as follows. 
The averaging block nonlinear Kaczmarz methods are proposed in Section \ref{secabnkm}. The convergence theory of the proposed methods with constant stepsize and adaptive stepsize is established in Section \ref{secconvanal}. 
Numerical experiments are provided to show the effectiveness of the proposed methods in Section \ref{secnumer}.
Finally, some remarks and conclusions are drawn in Section \ref{secconclu}.

\section{The averaging block nonlinear Kaczmarz methods}\label{secabnkm}
In this section, we first review some nonlinear Kaczmarz methods and then present a class of averaging block nonlinear Kaczmarz methods for the solution of the nonlinear system of equations \eqref{eqnesys}.

Consider the linear approximation of $F(x)$ based on the first-order Taylor expansion at $x_k$
\begin{equation*}
F(x)\approx F(x_k)+ F'(x_k)(x- x_k),
\end{equation*}
then the solution of nonlinear equations $F(x)=0$ can be approximated by a series of hyperplanes as follows
\begin{equation*}\label{eqhyplanes}
H_i= \{x\in \mathbb{R}^n: \langle \nabla F_i(x_k), x\rangle= -F_i(x_k)+ \langle \nabla F_i(x_k), x_k\rangle \},\quad i=1,\dots, m.
\end{equation*}
In the nonlinear Kaczmarz method \cite{2022WLB}, the approximate solution $x_{k+1}$ is obtained
by projecting the current vector $x_k$ onto the $i_k$-th hyperplane with the index $i_k$ is chosen from $[m]$.
 
When the index $i_k$ is selected with probability $p_{i_k}= {|F_{i_k}(x_k)|^2}/{\left\|F(x_k) \right\|^2_2}$, it gives the nonlinear randomized Kaczmarz method \cite{2022WLB}.
Assuming that there exists at least one vector $x_{\ast}$ such that $F(x_{\ast})=0$ and 
each nonlinear function $F_i(1\leq i\leq m)$ satisfies the local tangential cone condition in Definition \ref{defltcc}. Then, the convergence of the nonlinear randomized Kaczmarz method can be proved, which is restated in Theorem \ref{thmnrk}.
\begin{definition}{\em\cite{2007HLS} }\label{defltcc}
\rm For every $i\in\left\{1,2,\cdots, m \right\}$ and any $ x_1, x_2\in \mathbb{R}^n$, if there exists 
$\eta_i\in [0,\eta )(\eta= \max_i\eta_i< \frac{1}{2})$ such that
\begin{equation}
|F_i(x_1)- F_i(x_2)- (\nabla F_i(x_1)^T(x_1-x_2))| \leq \eta_i |F_i(x_1)- F_i(x_2) |,
\end{equation}
then the function $F:\mathbb{R}^n\rightarrow \mathbb{R}^m$ is referred to satisfy the local tangential cone condition.
\end{definition}
 
\begin{theorem}{\em\cite{2022WLB}}\label{thmnrk}
Consider that the nonlinear systems of equations $F(x)=0$ have a nonlinear function $F:D\subseteq\mathbb{R}^n\rightarrow \mathbb{R}^m$ on a bounded closed set $D$ and there  exists $x_{\ast}$ such that $F(x_{\ast})=0$. Assume that the derivative of $F$ is continuous in $D$ and $F'(x)$ is a row bounded below and full column rank matrix for any $x\in D$. If every nonlinear function $F_i ({1\leq i\leq m})$ satisfies the local tangential cone condition, $\eta= \max\limits_i\eta_i<\frac{1}{2}$, then the nonlinear randomized Kaczmarz method is convergent and it holds that
\begin{equation*}
    \mathbb{E}\left\| x_{k+1}-x_{\ast}\right\|_2^2\leq \left( 1-\frac{1-2\eta}{(1+\eta)^2m\kappa^2_F(F'(x_k))}\right)\mathbb{E}\left\| x_k-x_{\ast}\right\|_2^2,
\end{equation*}
where $\kappa_F(F'(x_k))=\left\|F'(x_k)\right\|_F /\sigma_{\min}(F'(x_k))$ is the scaled squared condition number of $F'(x_k)$ and $\sigma_{\min}(F'(x_k))$ is the smallest nonzero singular value of $F'(x_k)$. 
\end{theorem}
Denote $r_k= -F(x_k)$ be the residual vector with respect to $x_k$ and $r_k^{(i_k)}$ is its $i_k$-th component. By using a maximum residual control to choose the index $i_k$, that is,
$i_k= \arg \max\limits_{1\leq i\leq m} |r_k^{(i)}|^2$, a maximum residual nonlinear Kaczmarz method is presented and its convergence property is described in Theorem \ref{thmmrnk}.

\begin{theorem}{\em \cite{2023ZWZ} }\label{thmmrnk}
Under the same assumptions as Theorem \ref{thmnrk}, if every nonlinear function $F_i$ satisfies the local tangential cone condition, $\eta= \max\limits_i \eta_i<\frac{1}{2}$, then the maximum residual nonlinear Kaczmarz method is convergent and 
\begin{equation*}
\left\|x_{k+1}-x_{\ast} \right\|_2^2\leq \left(1-\frac{1-2\eta}{(1+\eta)^2m^2\kappa_F^2(F'(x_k))} \right) \left\|x_k-x_{\ast} \right\|_2^2.
\end{equation*}
\end{theorem}
Note that the above nonlinear Kaczmarz methods only use the information of a single nonlinear equation at each iteration, which will converge slowly when the size of the nonlinear system is very large.
To accelerate the convergence, it is natural to utilize multiple equations simultaneously at each iteration. Then, a block nonlinear Kaczmarz method \cite{2023ZWZ} is proposed and its iterate scheme is $x_{k+1}= x_{k}+ (F'_{\mathcal{I}_k}(x_k) )^{\dagger}F_{\mathcal{I}_k}(x_k)$ with $\mathcal{I}_k$ be the selected block index set. However, the calculation of the Moore-Penrose pesudoinverse of the submatrix $F'_{\mathcal{I}_k}(x_k)$ usually incurs a high cost.

In this work, to avoid computing the pesudoinverse of the submatrix, we introduce the averaging techniques which are also widely studied and used in the block Kaczmarz methods for linear systems, we refer the readers to \cite{2019N,2020DSS,2022MW,2023XYZ} for more details.
Specifically, for each selected block index set $\mathcal{I}_k\subset\{1,\cdots,m\}$, the block nonlinear Kaczmarz method with averaging technique iterates by
\begin{equation}\label{eqiterabnk}
 x_{k+1}= x_k- \alpha_k\left( \sum\limits_{i\in\mathcal{I}_k}\omega_k^{(i)} \frac{F_i(x_k)}{\left\|\nabla F_i(x_k) \right\| ^2_2} \nabla F_i(x_k) \right), \quad k \geq 0,
\end{equation}
where the weights satisfy $\omega_k^{(i)}\in[0,1]$ such that $\sum\limits_{i\in \mathcal{I}_k}\omega_k^{(i)} =1$ and $\alpha_k\in(0,2)$ is the stepsize.
Similar to \cite{2023ZWZ}, we consider the set $\mathcal{I}_k$ to be determined by the following greedy selection rule
\begin{equation}\label{eqabnkset}
 \mathcal{I}_k=\{ i_k| |r_k^{(i_k)}|^2\geq \theta\max\limits_{1\leq i\leq m} |r_k^{(i)}|^2 \}, \quad\theta\in(0,1].
\end{equation}

Then, based on the iterate scheme \eqref{eqiterabnk} and index subset selection criterion \eqref{eqabnkset}, we develop a class of averaging block nonlinear Kaczmarz methods, which is described in Algorithm \ref{algabnk}. 
\begin{algorithm}[h]
 \caption{Averaging block nonlinear Kaczmarz method} \label{algabnk}
 	\begin{algorithmic}[1]
  \Require $F, \ell, x_0$ and $\theta\in(0,1]$
    \Ensure $x_\ell$
        \State Initialize $x_0\in \mathbb{R}^n$ and $r_0= - F(x_0)$
 		\For {$k=0, 1,\ldots,\ell-1 $ }
   \State Determine the index set 
   $\mathcal{I}_k=\{i_k||r_k^{(i_k)}|^2\geq\theta\max\limits_{1\leq i\leq m}|r_k^{(i)}|^2\}$
         \State Update $x_{k+1}= x_k- \alpha_k\left( \sum\limits_{i\in\mathcal{I}_k}\omega_k^{(i)}\frac{F_i(x_k)}{\left\|\nabla F_i(x_k) \right\| ^2_2} \nabla F_i(x_k)  \right)$
        \State Compute residual $ r_k= -F(x_k)$
 		\EndFor
 	\end{algorithmic}
 \end{algorithm}

\begin{remark}\label{lemnonempty}
\rm The averaging block nonlinear Kaczmarz methods are well defined since the control index set $\mathcal{I}_k\subset\{1,2,\cdots, m\}$ in \eqref{eqiterabnk} selected by the approximate greedy criterion \eqref{eqabnkset} is nonempty. Indeed, if 
\begin{equation*}
     |r^{(i_s)}|^2 = \max\limits_{1\leq i\leq m} |r^{(i)}|^2, 
\end{equation*}
then $i_s\in\mathcal{I}_k$ for all $k\geq 0$. Thus, there must be an index in $\mathcal{I}_k$.
Moreover, the cardinality of the set $\mathcal{I}_k$ is flexible in the iteration process, that is, the component of the indices set $\mathcal{I}_k$ is changing with the increase of the iteration step $k$.
\end{remark}
Note that when the parameter $\theta=1$ in \eqref{eqabnkset} and $\alpha_k=1$ for all $k$, the averaging block nonlinear Kaczmarz method can recover the maximum residual nonlinear Kaczmarz method.
Moreover, the Algorithm \ref{algabnk} provides a general framework of averaging block nonlinear Kaczmarz methods for solving the nonlinear system of equations. 
By using suitable weights $\omega_k$ and stepsizes $\alpha_k$, it yields a series of averaging block nonlinear Kaczmarz methods.  
 
In Algorithm \ref{algabnk}, a simple choice for the weights is
$\omega_k^{(i)}=  {\left\|\nabla F_i(x_k) \right\|_2^2 }/\sum_{i\in\mathcal{I}_k} {\left\|\nabla F_i(x_k) \right\|_2^2}\in(0,1)$ for all $k\geq 0$, 
and then we can obtain the following update 
\begin{equation}\label{eqspupdate} 
  x_{k+1}= x_k - \alpha_k\frac{(F'_{\mathcal{I}_k}(x_k))^T F_{\mathcal{I}_k}(x_k)}{\left\|F'_{\mathcal{I}_k}(x_k) \right\|^2_2}. 
\end{equation} 
For the stepsize $\alpha_k$, we mainly focus on two choices:
\begin{itemize}
\item constant stepsize, that is, $\alpha_k$ is equal to a positive constant $\alpha \in(0,2)$ for all $k\geq 0$.
\item adaptive stepsize, which is based on the idea of extrapolation
\begin{equation}\label{eqadpsize}
\alpha_k= \delta \frac{\sum_{i\in\mathcal{I}_k}\hat{\omega}_k^{(i)}(F_i(x_k))^2}{\left\|\sum_{i\in\mathcal{I}_k}\hat{\omega}_k^{(i)}F_i(x_k)(\nabla F_i(x_k))^T \right\|_2^2},\quad k\geq 0,
\end{equation}
where $\hat{\omega}_k^{(i)}=\omega_k^{(i)}/\|\nabla F_i(x_k) \|_2^2$ and $\delta\in[1,2)$.
This extrapolated stepsize is also introduced for convex feasibility problems\cite{2019NRP} and for linear systems of equations\cite{2019N}.
\end{itemize}

\section{ Convergence analysis}\label{secconvanal}
In this section, the convergence properties of the averaging block nonlinear Kaczmarz methods with both constant stepsize and adaptive stepsize are discussed, respectively.
First, we introduce some lemmas that are crucial for the following convergence analysis of the proposed methods.
\begin{lemma}{\em \cite{2022ZLT}}\label{lemftau}
\rm If the function $F$ satisfies the local tangential cone condition, then for any $x_1, x_2\in \mathbb{R}^n$ and an index subset $\mathcal{I}\subset \left\{1,2,\cdots, m \right\}$, it holds that
\begin{equation}\label{lemftauneq1}
\left\| F_{\mathcal{I}}(x_1)- F_{\mathcal{I}}(x_2)- F'_{\mathcal{I}}(x_1)(x_1-x_2) \right\|^2_2 \leq \eta^2  \left\|F_{\mathcal{I}}(x_1)- F_{\mathcal{I}}(x_2)\right\|^2_2  
\end{equation}
and
\begin{equation}\label{lemftauneq2}
    \left\|F_{\mathcal{I}}(x_1)- F_{\mathcal{I}}(x_2)  \right\|_2^2\geq \frac{1}{1+\eta^2}\left\| F'_{\mathcal{I}}(x_1)(x_1- x_2) \right\|_2^2.
\end{equation}
\end{lemma}

\begin{lemma}\label{lemmatsv}
\rm Let $A\in \mathbb{R}^{m\times n}$ be any nonzero real matrix. For every vector $u\in  {\rm range}(A)$,  
\begin{equation*}
  \sigma_{\min}^2(A) \left\|u \right\|_2^2\leq \left\|A^T u \right\|^2_2 \leq \sigma_{\max}^2(A) \left\|u \right\|_2^2,
\end{equation*}
where ${\rm range}(A)$, $\sigma_{\min}(A)$  and $\sigma_{\max}(A)$ are the column space, the nonzero minimum and maximum singular values of $A$, respectively.
\end{lemma}

\subsection{Averaging block nonlinear Kaczmarz method with constant stepsize}
In this subsection, the averaging block nonlinear Kaczmarz method with constant stepsize $\alpha_k=\alpha$ and weights $\omega_k^{(i)}= {\left\|\nabla F_i(x_k) \right\|_2^2 }/\sum_{i\in\mathcal{I}_k} {\left\|\nabla F_i(x_k) \right\|_2^2}$ is considered and its convergence analysis is given. 
In this case, the update formula of the averaging block nonlinear Kaczmarz method in Algorithm \ref{algabnk} becomes 
\begin{equation}\label{eqiterabnkcon}
   x^{k+1}= x^{k}- \alpha \frac{(F'_{\mathcal{I}_k}(x_k))^T F_{\mathcal{I}_k}(x_k) }{\left\|F'_{\mathcal{I}_k}(x_k) \right\|^2_2 }. 
\end{equation}

\begin{lemma}\label{lemcstabnkiter}
 \rm  If the function $F$ satisfies the local tangential cone condition, $\eta=\max\limits_i \eta_i< \frac{1}{2}$ and there exists a vector $x_{\ast}$ satisfies $F(x_{\ast})=0$,  then from the iteration formula \eqref{eqiterabnkcon}, $\alpha\in(0,2(1-\eta))$ and $\mathcal{I}_k\subset \left\{1,2,\cdots, m \right\} $, it satisfies that
\begin{equation} 
\left\|x_{k+1}-x_{\ast} \right\|_2^2\leq \left\|x_{k}-x_{\ast} \right\|_2^2-(2(1-\eta)\alpha-\alpha^2)\frac{\left\|F_{\mathcal{I}_k}(x_k)\right\|_2^2}{\left\| F'_{\mathcal{I}_k}(x_k)\right\|_2^2}.
\end{equation}
\end{lemma}

\begin{Proof}
From the iteration \eqref{eqiterabnkcon}, it follows that
\begin{equation*}
\begin{aligned}
&\left\|x_{k+1}-x_{\ast} \right\|_2^2=\left\|x_{k}-x_{\ast} \right\|_2^2
 +\left\|x_{k+1}-x_{k} \right\|_2^2+ 2\langle x_{k+1}-x_{k}, x_{k}-x_{\ast} \rangle\\
&=\left\|x_{k}-x_{\ast} \right\|_2^2+\left\|\alpha\frac{(F'_{\mathcal{I}_k}(x_k))^T F_{\mathcal{I}_k}(x_k)}{\left\|F'_{\mathcal{I}_k}(x_k) \right\|^2_2 }\right\|^2_2 
+ 2\left \langle -\alpha \frac{(F'_{\mathcal{I}_k}(x_k))^TF_{\mathcal{I}_k}(x_k)}{\left\|F'_{\mathcal{I}_k}(x_k)\right\|^2_2},  x_{k}-x_{\ast} \right \rangle\\
 &=\left\|x_{k}-x_{\ast} \right\|_2^2+ \alpha^2\left\|\frac{(F'_{\mathcal{I}_k}(x_k))^TF_{\mathcal{I}_k}(x_k)}{\left\|F'_{\mathcal{I}_k}(x_k)\right\|_2^2} \right\|_2^2
- 2\alpha \frac{  F^T_{\mathcal{I}_k}(x_k)F'_{\mathcal{I}_k}(x_k)}{\left\|F'_{\mathcal{I}_k}(x_k)\right\|_2^2}(x_k-x_{\ast})\\
 &=\left\|x_{k}-x_{\ast} \right\|_2^2+ (\alpha^2-2\alpha)\frac{\left\|F_{\mathcal{I}_k}(x_k)\right\|_2^2}{\left\| F'_{\mathcal{I}_k}(x_k) \right\|_2^2} +2\alpha\frac{F^T_{\mathcal{I}_k}(x_k)}{\left\| F'_{\mathcal{I}_k}(x_k) \right\|_2^2}\left( F_{\mathcal{I}_k}(x_k)-F_{\mathcal{I}_k}(x_{\ast})-F'_{\mathcal{I}_k}(x_k)(x_k-x_{\ast}) \right).
\end{aligned}
\end{equation*}
Moreover, using the Cauchy-Schwarz inequality, it obtains that
\begin{equation*}
\begin{aligned}
 \left\|x_{k+1}-x_{\ast} \right\|_2^2=&\left\|x_{k}-x_{\ast} \right\|_2^2+(\alpha^2-2\alpha)\frac{\left\|F_{\mathcal{I}_k}(x_k)\right\|_2^2}{\left\| F'_{\mathcal{I}_k}(x_k)\right\|_2^2}\\
&+ 2\alpha\frac{\left\|F^T_{\mathcal{I}_k}(x_k)\right\|_2}{\left\| F'_{\mathcal{I}_k}(x_k) \right\|_2^2}\left\|F_{\mathcal{I}_k}(x_k)-F_{\mathcal{I}_k}(x_{\ast})-F'_{\mathcal{I}_k}(x_k)(x_k-x_{\ast})\right\|_2.
\end{aligned}
\end{equation*}
From \eqref{lemftauneq1} in Lemma \ref{lemftau}, it yields that
\begin{equation}
\begin{aligned}
 \left\|x_{k+1}-x_{\ast} \right\|_2^2
 &\leq \left\|x_{k}-x_{\ast} \right\|_2^2 
 +(\alpha^2-2\alpha)\frac{\left\|F_{\mathcal{I}_k}(x_k)\right\|_2^2}{\left\| F'_{\mathcal{I}_k}(x_k)\right\|_2^2}+ 2\alpha\eta \frac{\left\|F^T_{\mathcal{I}_k}(x_k)\right\|_2}{\left\| F'_{\mathcal{I}_k}(x_k) \right\|_2^2}\left\|F_{\mathcal{I}_k}(x_k)-F_{\mathcal{I}_k}(x_{\ast}) \right\|_2 \\
 &\leq\left\|x_{k}-x_{\ast} \right\|_2^2 
 +(\alpha^2-2\alpha)\frac{\left\|F_{\mathcal{I}_k}(x_k)\right\|_2^2}{\left\| F'_{\mathcal{I}_k}(x_k)\right\|_2^2}+ 2\alpha\eta \frac{\left\|F^T_{\mathcal{I}_k}(x_k)\right\|_2}{\left\| F'_{\mathcal{I}_k}(x_k) \right\|_2^2}\left\|F_{\mathcal{I}_k}(x_k)\right\|_2 \\
 &\leq\left\|x_{k}-x_{\ast} \right\|_2^2 
 -(2(1-\eta)\alpha-\alpha^2)\frac{\left\|F_{\mathcal{I}_k}(x_k)\right\|_2^2}{\left\| F'_{\mathcal{I}_k}(x_k)\right\|_2^2},
\end{aligned}
\end{equation}
where the second inequality depends on $F(x_{\ast})=0$.

\end{Proof}
The Lemma \ref{lemcstabnkiter} provides an important property for the sequence $\{x_k\}_{k=0}^{\infty}$ generated by the scheme \eqref{eqiterabnkcon} at each iteration. Next, we give an upper bound of the convergence rate for the averaging block nonlinear Kaczmarz method with constant stepsize, which is described in Theorem \ref{thmcstABNK}.

\begin{theorem}\label{thmcstABNK}
Under the same condition as Theorem \ref{thmnrk}, for every $i\in\left\{1,\cdots, m \right\}$, if the nonlinear function $F_i$ satisfies the local tangential cone condition, $\eta= \max\limits_i \eta_i<\frac{1}{2}$ and exists a vector such that $F(x_{\ast})=0$, then the averaging block nonlinear Kaczmarz method with constant stepsize $\alpha\in[1, 2(1-\eta))$ is convergent and 
\begin{equation}
    \left\|x_{k+1}-x_{\ast} \right\|_2^2\leq  \left(1- \frac{ 2(1-\eta)\alpha-\alpha^2 }{(1+\eta^2)\kappa_F^2(F'_{\mathcal{I}_k}(x_k))} \right)\left\|x_{k}-x_{\ast} \right\|_2^2.
\end{equation}
\end{theorem}

\begin{Proof}
From Lemma \ref{lemcstabnkiter}, it follows that
\begin{equation*}
 \begin{aligned}
\left\|x_{k+1}-x_{\ast} \right\|_2^2
&\leq \left\|x_{k}-x_{\ast} \right\|_2^2- (2(1-\eta)\alpha-\alpha^2)\frac{\left\|F_{\mathcal{I}_k}(x_k)-F_{\mathcal{I}_k}(x_{\ast})\right\|_2^2}{\left\| F'_{\mathcal{I}_k}(x_k)\right\|_2^2}\\
&\leq \left\|x_{k}-x_{\ast} \right\|_2^2- \frac{2(1-\eta)\alpha-\alpha^2}{(1+\eta^2)\left\| F'_{\mathcal{I}_k}(x_k)\right\|_2^2}\left\| F'_{\mathcal{I}_k}(x_k)(x_k- x_{\ast}) \right\|_2^2\\
& \leq \left\|x_{k}-x_{\ast} \right\|_2^2- \frac{2(1-\eta)\alpha-\alpha^2}{(1+\eta^2)\left\| F'_{\mathcal{I}_k}(x_k)\right\|_2^2} \sigma_{\min}^2(F'_{\mathcal{I}_k}(x_k))\left\|x_{k}-x_{\ast} \right\|_2^2\\ 
& \leq \left\|x_{k}-x_{\ast} \right\|_2^2- \frac{2(1-\eta)\alpha-\alpha^2}{(1+\eta^2)\left\| F'_{\mathcal{I}_k}(x_k)\right\|_F^2} \sigma_{\min}^2(F'_{\mathcal{I}_k}(x_k))\left\|x_{k}-x_{\ast} \right\|_2^2 \\ 
&\leq \left\|x_{k}-x_{\ast} \right\|_2^2-\frac{2(1-\eta)\alpha-\alpha^2 }{(1+\eta^2)\kappa_F^2(F'_{\mathcal{I}_k}(x_k))}\left\|x_{k}-x_{\ast} \right\|_2^2\\
& \leq \left(1- \frac{2(1-\eta)\alpha-\alpha^2}{(1+\eta^2)\kappa_F^2(F'_{\mathcal{I}_k}(x_k))} \right)\left\|x_{k}-x_{\ast} \right\|_2^2.
\end{aligned} 
\end{equation*}
The second and third inequality hold because \eqref{lemftauneq2} and Lemma \eqref{lemmatsv}, respectively. 
\end{Proof}

\begin{remark}
\rm From Theorem \ref{thmcstABNK}, it is observed that the convergence factor of the averaging block nonlinear Kaczmarz method with constant stepsize is
\begin{equation*}
\rho(\alpha)= \  1- \frac{2(1-\eta)\alpha-\alpha^2}{(1+\eta^2)\kappa_F^2(F'_{\mathcal{I}_k}(x_k))},
\end{equation*}
which achieves its minimum $1-\frac{1-2\eta+\eta^2}{(1+\eta^2)\kappa_F^2(F'_{\mathcal{I}_k}(x_k))}$ at $\alpha= 1-\eta$ and is equal to $1-\frac{1-2\eta }{(1+\eta^2)\kappa_F^2(F'_{\mathcal{I}_k}(x_k))}$ when $\alpha=1$. 
\end{remark}

\subsection{Averaging block nonlinear Kaczmarz method with adaptive stepsize}
In this subsection, we analyze the convergence of the averaging block nonlinear Kaczmarz method with the adaptive stepsize $\alpha_k= \delta {\sum_{i\in\mathcal{I}_k}\hat{\omega}_k^{(i)}(F_i(x_k))^2}/{\left\|\sum_{i\in\mathcal{I}_k}\hat{\omega}_k^{(i)}F_i(x_k)(\nabla F_i(x_k))^T \right\|_2^2}$ in \eqref{eqadpsize} and weights $\omega_k^{(i)}= {\left\|\nabla F_i(x_k) \right\|_2^2 }/\sum_{i\in\mathcal{I}_k} {\left\|\nabla F_i(x_k) \right\|_2^2}$. 

\begin{lemma}\label{lemadpabnkiter}
 \rm  If the function $F$ satisfies the local tangential cone condition and a vector $x_{\ast}$ satisfies $F(x_{\ast})=0$, $\eta=\max\limits_i \eta_i< \frac{1}{2}$, then from the averaging block iteration formula \eqref{eqiterabnk} with adaptive stepsize 
\eqref{eqadpsize}, $\delta\in (0,2(1-\eta))$ and $\mathcal{I}_k\subset \{1,2,\cdots,m\}$, it obeys that
    \begin{equation}\label{eqestiadp}
        \left\|x_{k+1}-x_{\ast} \right\|_2^2\leq \left\|x_{k}-x_{\ast} \right\|_2^2- (2(1-\eta)\delta-\delta^2 ) \frac{\left\|F_{\mathcal{I}_k}(x_k) \right\|_2^2}{\sigma^2_{\max}(F'_{\mathcal{I}_k}(x_k))}.
    \end{equation} 
\begin{Proof}
It follows that
\begin{equation*}
\begin{aligned}
&\left\|x_{k+1}-x_{\ast} \right\|_2^2
=\left\|x_{k}-x_{\ast} \right\|_2^2 + \left\|x_{k+1}-x_{k} \right\|_2^2+ 2\langle x_{k+1}-x_{k}, x_{k}-x_{\ast} \rangle \\
&=\left\|x_{k}-x_{\ast} \right\|_2^2 + \left\|\alpha_k\frac{(F'_{\mathcal{I}_k}(x_k))^TF_{\mathcal{I}_k}(x_k)}{\left\|F'_{\mathcal{I}_k}(x_k) \right\|^2_2 }\right\|^2_2 
+ 2\left \langle -\alpha_k\frac{(F'_{\mathcal{I}_k}(x_k))^TF_{\mathcal{I}_k}(x_k)}{\left\|F'_{\mathcal{I}_k}(x_k)\right\|^2_2}, x_{k}-x_{\ast} \right \rangle\\
&=\left\|x_{k}-x_{\ast} \right\|_2^2 + \left\|\delta \frac{\sum_{i\in\mathcal{I}_k}\hat{\omega}_k^{(i)}(F_i(x_k))^2}{\left\|\sum_{i\in\mathcal{I}_k}\hat{\omega}_k^{(i)}F_i(x_k)(\nabla F_i(x_k))^T \right\|_2^2}\frac{(F'_{\mathcal{I}_k}(x_k))^TF_{\mathcal{I}_k}(x_k)}{\left\|F'_{\mathcal{I}_k}(x_k) \right\|^2_2 }\right\|^2_2 \\
&\quad - 2\left \langle \delta \frac{\sum_{i\in\mathcal{I}_k}\hat{\omega}_k^{(i)}(F_i(x_k))^2}{\left\|\sum_{i\in\mathcal{I}_k}\hat{\omega}_k^{(i)}F_i(x_k)(\nabla F_i(x_k))^T \right\|_2^2}\frac{(F'_{\mathcal{I}_k}(x_k))^TF_{\mathcal{I}_k}(x_k)}{\left\|F'_{\mathcal{I}_k}(x_k)\right\|^2_2},  x_{k}-x_{\ast} \right \rangle \\
&=\left\|x_{k}-x_{\ast} \right\|_2^2 + \delta^2\left\|\frac{\left\|F_{\mathcal{I}_k}(x_k)\right\|_2^2(F'_{\mathcal{I}_k}(x_k))^TF_{\mathcal{I}_k}(x_k)}{\left\|(F'_{\mathcal{I}_k}(x_k) )^TF_{\mathcal{I}_k}(x_k)\right\|_2^2} \right\|_2^2
- 2\delta \frac{ \left\|F_{\mathcal{I}_k}(x_k)\right\|_2^2F^T_{\mathcal{I}_k}(x_k)F'_{\mathcal{I}_k}(x_k)}{\left\|(F'_{\mathcal{I}_k}(x_k) )^TF_{\mathcal{I}_k}(x_k)\right\|_2^2}(x_k-x_{\ast}) \\
&=\left\|x_{k}-x_{\ast} \right\|_2^2 +2\delta\frac{\left\|F_{\mathcal{I}_k}(x_k)\right\|_2^2F^T_{\mathcal{I}_k}(x_k)}{\left\|(F'_{\mathcal{I}_k}(x_k) )^Tr_k^{(\mathcal{I}_k)}\right\|_2^2}( F_{\mathcal{I}_k}(x_k)-F_{\mathcal{I}_k}(x_{\ast})-F'_{\mathcal{I}_k}(x_k)(x_k-x_{\ast}) )\\
 &\quad - 2\delta\frac{\left\|F_{\mathcal{I}_k}(x_k)\right\|_2^2F^T_{\mathcal{I}_k}(x_k)}{\left\|(F'_{\mathcal{I}_k}(x_k) )^TF_{\mathcal{I}_k}(x_k)\right\|_2^2} F_{\mathcal{I}_k}(x_k)+  \delta^2\frac{\left\|r_k^{(\mathcal{I}_k)}\right\|_2^4}{\left\|(F'_{\mathcal{I}_k}(x_k) )^Tr_k^{(\mathcal{I}_k)}\right\|_2^2} \\
&=\left\|x_{k}-x_{\ast} \right\|_2^2 +2\delta\frac{\left\|F_{\mathcal{I}_k}(x_k)\right\|_2^2F^T_{\mathcal{I}_k}(x_k)}{\left\|(F'_{\mathcal{I}_k}(x_k) )^TF_{\mathcal{I}_k}(x_k)\right\|_2^2}\left( F_{\mathcal{I}_k}(x_k)-F_{\mathcal{I}_k}(x_{\ast})-F'_{\mathcal{I}_k}(x_k)(x_k-x_{\ast}) \right)\\
&\quad + (\delta^2-2\delta)\frac{\left\|F_{\mathcal{I}_k}(x_k)\right\|_2^4}{\left\|(F'_{\mathcal{I}_k}(x_k) )^TF_{\mathcal{I}_k}(x_k)\right\|_2^2}.
\end{aligned}
\end{equation*}
Furthermore, using the Cauchy-Schwarz inequality, it holds that
\begin{equation*}
\begin{aligned}
 \left\|x_{k+1}-x_{\ast} \right\|_2^2&\leq \left\|x_{k}-x_{\ast} \right\|_2^2
  +2\delta\frac{\left\|F_{\mathcal{I}_k}(x_k) \right\|_2^2\left\|F^T_{\mathcal{I}_k}(x_k)\right\|_2}{\left\|(F'_{\mathcal{I}_k}(x_k) )^TF_{\mathcal{I}_k}(x_k)\right\|_2^2}\left\|F_{\mathcal{I}_k}(x_k)-F_{\mathcal{I}_k}(x_{\ast})-F'_{\mathcal{I}_k}(x_k)(x_k-x_{\ast})\right\|_2\\
&\quad + (\delta^2-2\delta)\frac{\left\|F_{\mathcal{I}_k}(x_k)\right\|_2^4}{\left\|(F'_{\mathcal{I}_k}(x_k) )^TF_{\mathcal{I}_k}(x_k)\right\|_2^2}.
    \end{aligned}
\end{equation*}
From \eqref{lemftauneq1} in Lemma \ref{lemftau} and $F(x_{\ast})=0$, it yields that
\begin{equation}
\begin{aligned}
 \left\|x_{k+1}-x_{\ast} \right\|_2^2
 &\leq \left\|x_{k}-x_{\ast} \right\|_2^2 
 + 2\delta\eta \frac{\left\|F_{\mathcal{I}_k}(x_k) \right\|_2^2}{\left\|(F'_{\mathcal{I}_k}(x_k) )^T F_{\mathcal{I}_k}(x_k)\right\|_2^2}\left\|F^T_{\mathcal{I}_k}(x_k)\right\|_2 \left\|F_{\mathcal{I}_k}(x_k) \right\|_2\\
 &\quad + (\delta^2-2\delta)\frac{\left\|F_{\mathcal{I}_k}(x_k)\right\|_2^4}{\left\|(F'_{\mathcal{I}_k}(x_k) )^TF_{\mathcal{I}_k}(x_k) \right\|_2^2}\\
 &\leq \left\|x_{k}-x_{\ast} \right\|_2^2- \left(2(1-\eta)\delta-\delta^2\right) \frac{\left\|F_{\mathcal{I}_k}(x_k)\right\|_2^4}{\left\|(F'_{\mathcal{I}_k}(x_k) )^TF_{\mathcal{I}_k}(x_k)\right\|_2^2}.
\end{aligned}
\end{equation}
From Lemma \ref{lemmatsv}, it has
\begin{equation}
    \begin{aligned}
  \left\|x_{k+1}-x_{\ast} \right\|_2^2
  &\leq \left\|x_{k}-x_{\ast} \right\|_2^2- \left(2(1-\eta)\delta-\delta^2\right) \frac{\left\|F_{\mathcal{I}_k}(x_k)\right\|_2^4}{\sigma_{\max}^2(F'_{\mathcal{I}_k}(x_k) )\left\| F_{\mathcal{I}_k}(x_k)\right\|_2^2}\\
  &\leq \left\|x_{k}-x_{\ast} \right\|_2^2-\left(2(1-\eta)\delta-\delta^2\right) \frac{\left\|F_{\mathcal{I}_k}(x_k)\right\|_2^2}{\sigma_{\max}^2(F'_{\mathcal{I}_k}(x_k) )},
    \end{aligned}
\end{equation}
which is exactly the estimate in \eqref{eqestiadp}.
\end{Proof}
\end{lemma}
Lemma \ref{lemadpabnkiter} shows a crucial relationship between iterate sequence $\{x_k\}_{k=0}^{\infty}$ and solution vector $x_{\ast}$. In Theorem \ref{thmadpABNK}, we give an upper bound of the convergence rate for the averaging block nonlinear Kaczmarz method with adaptive stepsize defined by \eqref{eqadpsize}.

\begin{theorem}\label{thmadpABNK}
Under the same condition as Theorem \ref{thmnrk}, for every $i\in\{1,\cdots,m\}$, if the nonlinear function $F_i$ satisfies the local tangential cone condition, $\eta= \max\limits_i \eta_i<\frac{1}{2}$ and exists a vector such that $F(x_{\ast})=0$, then the averaging block nonlinear Kaczmarz method with adaptive stepsize and $\delta\in(0, 2(1-\eta))$ is convergent. Moreover, it satisfies that 
\begin{equation}
    \left\|x_{k+1}-x_{\ast} \right\|_2^2\leq \left(1- \frac{2(1-\eta)\delta-\delta^2 }{(1+\eta^2)\kappa^2(F'_{\mathcal{I}_k}(x_k)) } \right) \left\|x_{k}-x_{\ast} \right\|_2^2.
\end{equation} 
\end{theorem}

\begin{Proof}
From Lemma \ref{lemadpabnkiter}, it follows that
\begin{equation}
 \begin{aligned}
\left\|x_{k+1}-x_{\ast} \right\|_2^2
&\leq \left\|x_{k}-x_{\ast} \right\|_2^2- \frac{2(1-\eta)\delta-\delta^2}{\sigma^2_{\max}(F'_{\mathcal{I}_k}(x_k))} \left\|F_{\mathcal{I}_k}(x_k)-F_{\mathcal{I}_k}(x_{\ast}) \right\|_2^2 \\
&\leq \left\|x_{k}-x_{\ast} \right\|_2^2-\frac{2(1-\eta)\delta-\delta^2}{(1+\eta^2)\sigma^2_{\max}(F'_{\mathcal{I}_k}(x_k)) }\left\|F'_{\mathcal{I}_k}(x_k)(x_k- x_{\ast}) \right\|_2^2\\
&\leq \left\|x_{k}-x_{\ast} \right\|_2^2-\frac{2(1-\eta)\delta-\delta^2}{(1+\eta^2)\sigma^2_{\max}(F'_{\mathcal{I}_k}(x_k)) }\sigma_{\min}^2(F'_{\mathcal{I}_k}(x_k))\left\|x_{k}-x_{\ast} \right\|_2^2\\
& \leq \left(1- \frac{2(1-\eta)\delta-\delta^2}{(1+\eta^2)\kappa^2(F'_{\mathcal{I}_k}(x_k)) } \right)\left\|x_{k}-x_{\ast} \right\|_2^2.
\end{aligned} 
\end{equation}
The second inequality holds because \eqref{lemftauneq2} and the last inequality depends on the definition of condition number $\kappa(\cdot)$. 
\end{Proof} 

\begin{remark}
\rm From Theorem \ref{thmadpABNK}, it is observed that the convergence factor of the averaging block nonlinear Kaczmarz method with adaptive stepsize is
\begin{equation*}
\rho(\delta)= \  1- \frac{2(1-\eta)\delta-\delta^2}{(1+\eta^2)\kappa^2(F'_{\mathcal{I}_k}(x_k))},
\end{equation*}
which achieves its minimum $1-\frac{1-2\eta+\eta^2}{(1+\eta^2)\kappa^2(F'_{\mathcal{I}_k}(x_k))}$ at $\delta= 1-\eta$ and is equal to $1-\frac{1-2\eta }{(1+\eta^2)\kappa^2(F'_{\mathcal{I}_k}(x_k))}$ when $\delta=1$. 
\end{remark}

\begin{remark}\label{remabnkconver}
\rm Moreover, Theorems \ref{thmnrk} and \ref{thmmrnk} show that the convergence factors of the nonlinear randomized Kaczmarz method and maximum residual nonlinear Kaczmarz method are 
$$1-\frac{1-2\eta}{(1+\eta)^2m^2\kappa_F^2(F'(x_k))}
\text{ and }  1-\frac{1-2\eta} {(1+\eta)^2m\kappa^2_F(F'(x_k))},$$
respectively.

Therefore, from Theorems \ref{thmcstABNK}, \ref{thmadpABNK} and the inequality $\kappa^2(A)\leq \kappa^2_F(A)$, the theoretical upper bound of the convergence rate for the averaging block nonlinear Kaczmarz method is larger than those of the nonlinear randomized Kaczmarz and the maximum residual nonlinear Kaczmarz methods if the parameters are chosen appropriately. 
\end{remark}

\section{Numerical experiments} \label{secnumer}
 In this section, numerical experiments are presented to verify the efficiency of averaging block nonlinear Kaczmarz methods from the aspects of the number of iteration steps (denoted as `IT') and the elapsed CPU time in seconds (denoted as `CPU'). 
 
 In the experiments, the averaging block nonlinear Kaczmarz methods with constant stepsize (abbreviated as `ABNK-1') and adaptive stepsize (abbreviated as `ABNK-2') are compared with the nonlinear randomized Kaczmarz method (abbreviated as `NRK'), the maximum residual nonlinear Kaczmarz method (abbreviated as `MRNK') and the maximum residual block nonlinear Kaczmarz method (abbreviated as `MRBNK').
  
In the maximum residual block nonlinear Kaczmarz method, the LSQR algorithm \cite{1982PS} is used to solve the least squares problem instead of directly calculating the Moore-Penrose pesudoinverse at each iteration. 
The optimal parameters for the maximum residual block Kaczmarz method and averaging block nonlinear Kaczmarz methods are experimentally selected by minimizing the number of iteration steps.
 
All the iterations terminate when the residual $\left\|r_k \right\|^2_2$ (denoted as `RES') satisfies RES$\leq 10^{-6}$ or the number of iterations exceeds a maximum number, e.g., 400,000.
 
\subsection{H-equation problem \cite{2010BDDKB}}\label{exam1}
\begin{equation}\label{eqHequation}
F(H)(\mu)= H(\mu)- \left(1- \frac{c}{2}\int_0^1 \frac{\mu H(\nu)d \nu }{\mu + \nu} \right)^{-1} =0,
\end{equation}
which arises from radiative transfer theory. When the integrals is approximated by the composite midpoint rule
\begin{equation}\label{eqmidrule}
\int_0^1 f(\mu) d\mu\approx \frac{1}{n} \sum\limits_{j=1}^{n} f(\mu_j),
\end{equation}
where $\mu_j= (j-\frac{1}{2})/n, (j=1,2,\cdots, n)$, the following discrete problem is obtained.
\begin{equation}\label{eqHeqdiscrete}
F(x)_i= x_i- \left(1- \frac{c}{2n} \sum\limits_{j=1}^{n} \frac{\mu_ix_j}{\mu_i+\mu_j}  \right)^{-1}, \quad i=1,2,\cdots, n.
\end{equation}
In this experiment, the initial vector is $x_0=(0,0,\cdots, 0)^T$ and the constant $c= 0.9$. 
 
In Table \ref{tab:resultEx1}, the number of iteration steps and the elapsed CPU time for all tested nonlinear Kaczmarz-type methods are listed when the number of equations $m$ varies, respectively.

From Table \ref{tab:resultEx1}, it is observed that all methods can successfully converge to the solution of the nonlinear problem. 
Both the number of iteration steps and elapsed CPU time of the ABNK-1 and ABNK-2 methods increase slowly as the size of the problem increases, which indicates that the ABNK methods are advantageous for solving large nonlinear systems of equations.
It is also seen that the ABNK-2 method requires the fewest number of iteration steps and least elapsed CPU time among all methods. The number of iteration steps and the CPU time of ABNK-1 method are similar to those of the MRBNK method. 
These results imply that the average scheme is effective and can greatly improve the convergence speed of the nonlinear Kaczmarz method when suitable stepsizes are given. 

Moreover, it is computed from Table \ref{tab:resultEx1} that the speed-up of the ABNK-2 method against the NRK method is around 10, the speed-up of the ABNK-2 method against the MRNK method is around 8, the speed-up of the ABNK-2 method against the MRBNK method nears to 3 and the speed-up of the ABNK-2 method against the ABNK-1 method is distributed nears to 3. It implies that the weights and the adaptive stepsize selected in the ABNK-2 method are more efficient than those of other methods.

\begin{table} [!htbp] 
\centering 
\caption{ Numerical results for Example \ref{exam1}. }\label{tab:resultEx1}  
\resizebox{\textwidth}{!}{  
\begin{tabular}{lcccccccccc} 
\hline
\multirow{2}{*} {$m$}  & \multicolumn{2}{c}{NRK} & \multicolumn{2}{c}{MRNK} & \multicolumn{2}{c}{MRBNK} & \multicolumn{2}{c}{ABNK-1} & \multicolumn{2}{c}{ABNK-2} \\  
 \cmidrule[0.25mm](lr){2-3} \cmidrule[0.25mm](lr){4-5} \cmidrule[0.25mm](lr){6-7}\cmidrule[0.25mm](lr){8-9} \cmidrule[0.25mm](lr){10-11}  & IT &CPU &IT &CPU &IT &CPU &IT &CPU &IT &CPU  \\ 
\hline
100 &	 2017&	 0.268&	 1808&	 0.217&	 21&	 0.078&	20&	 0.073 &	  12&	 0.027  \\  
200 &	 4238&	 1.974&	 3783&	 1.696&	 22&	 0.552&	22&	 0.584 &	  13&	 0.213  \\  
300 &	 6523&	 6.676&	 5820&	 5.778&	 22&	 1.716&	22&	 1.788 &	  13&	 0.686  \\  
400 &	 8859&	 15.833&	 7888&	 13.589&	 23&	 4.060&	22&	 4.105 &	  14&	 1.693  \\  
500 &	 11234&	 30.361&	 9998&	 26.793&	 23&	 7.826&	23&	 8.398 &	  14&	 3.264  \\  
600 &	 13608&	 53.229&	 12126&	 46.887&	 24&	 14.373&	23&	 14.325 &	  14&	 5.765  \\  
700 &	 16040&	 85.367&	 14273&	 76.062&	 24&	 23.086&	24&	 23.244 &	  14&	 9.154  \\  
800 &	 18431&	 127.871&	 16430&	 114.740&	 24&	 33.687&	24&	 34.618 &	  14&	 13.305  \\  
900 &	 20911&	 186.210&	 18594&	 161.518&	 24&	 48.406&	24&	 49.120 &	  14&	 19.383  \\  
1000 &	 23358&	 250.946&	 20786&	 222.507&	 24&	 66.670&	24&	 67.435 &	  14&	 26.259  \\  
\hline
\end{tabular}  
 }   
 \begin{tablenotes} 
 \footnotesize               
\item The parameters in the MRBNK, ABNK-1 and ABNK-2 methods are $\theta_1^*=0.1, (\alpha^*, \theta_2^*)=(1.7,0.1)$ and $(\delta^*, \theta_3^*)=(1.2,0.2)$, respectively. 
\end{tablenotes}
\end{table} 

\begin{figure}[!htbp] 
	\centering
	\subfigure[$m=100$]
	{
		\begin{minipage}[t]{0.48\linewidth}
			\centering
			\includegraphics[width=1\textwidth]{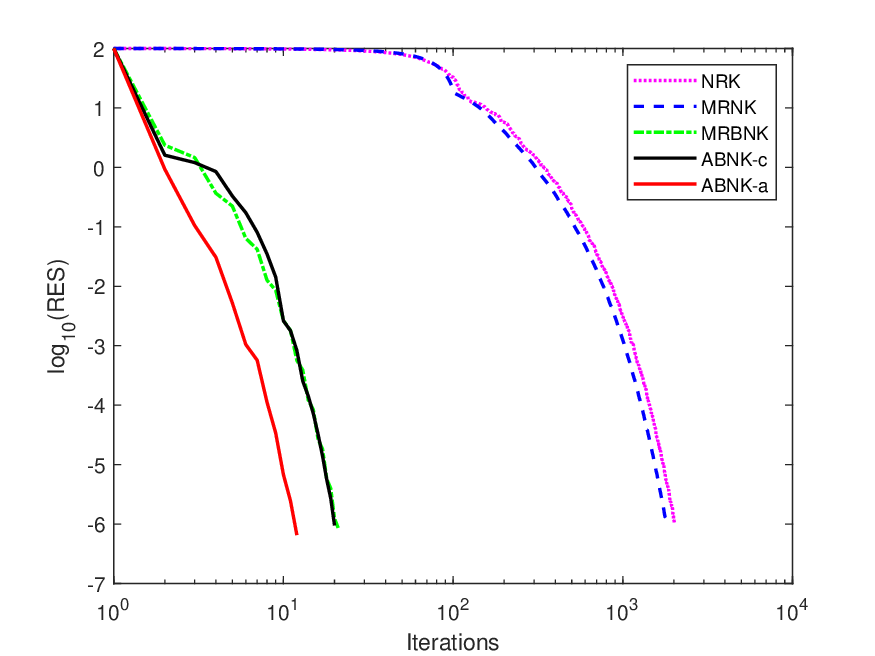}
		\end{minipage}
	}
	\subfigure[$m=100$]
	{
		\begin{minipage}[t]{0.48\linewidth}
			\centering
			\includegraphics[width=1\textwidth]{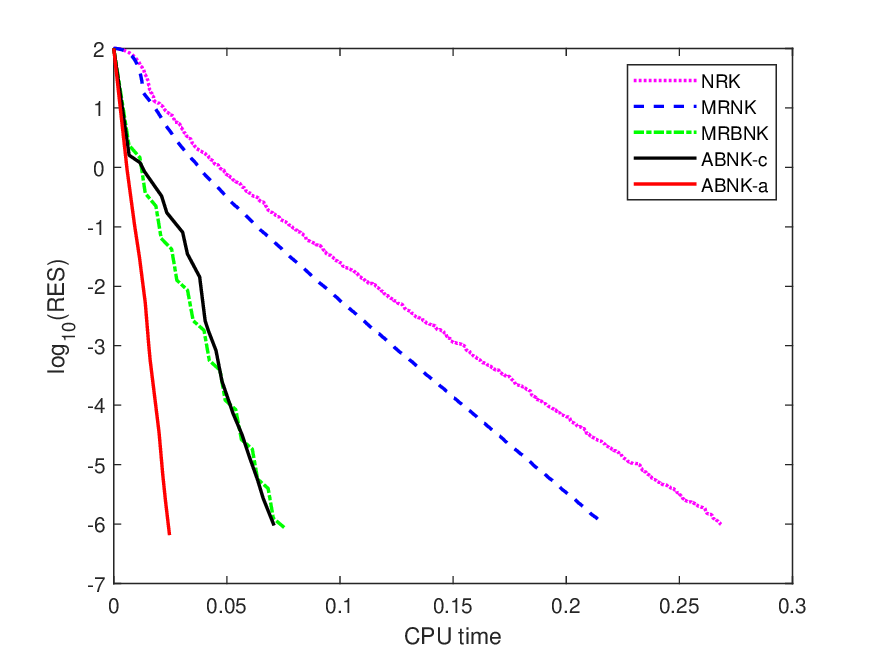}
		\end{minipage}
	}
	\caption{Curves of residual versus the iteration steps and CPU time in Example \ref{exam1}.} 	
		\label{fig:resvsitcpuEx1}
\end{figure}

In Figure \ref{fig:resvsitcpuEx1}, the curves of the residuals versus the number of iteration steps and the elapsed CPU time for all tested methods are depicted when $m=100$, respectively. 

From Figure \ref{fig:resvsitcpuEx1}, it is observed that the curves of residuals versus both the number of iteration steps and the elapsed CPU time for the ABNK methods decrease faster than those of other methods. Moreover, the curves for ABNK-2 method decrease the fastest among all methods.

\subsection{Tridiagonal system problem \cite{1989L}}\label{exam2}
\begin{equation}\label{eqtridiagonalsys}
\begin{cases}
F_k(x) = 4(x_k- x_{k+1}^2),   & k= 1,\\
F_k(x) = 8x_k(x_k^2-x_{k-1}) -2(1-x_k)+ 4(x_k- x_{k+1}^2),   & 1<k<n, \\
F_k(x) = 8x_k(x_k^2-x_{k-1}) -2(1-x_k),  & k=n, \\
m = n.
\end{cases}
\end{equation}
Here, $m$ and $n$ are the number of nonlinear equations and the dimension of vector $x$, respectively. In this experiment, the initial vector is $x_0= (12,12,\cdots,12)^T$. 
 
In Table \ref{tab:resultEx2}, the number of iteration steps and the elapsed CPU time for all  nonlinear Kaczmarz methods are listed when the size $m$ varies, respectively.

From Table \ref{tab:resultEx2}, it is observed that all methods can compute the solution of the nonlinear problem successfully.
The ABNK methods are faster than other three methods both in terms of the number of iteration steps and the CPU time. Moreover, the ABNK-2 method requires the fewest number of iteration steps and least elapsed CPU time among all methods. It further shows the efficiency of the adaptive stepsize chosen in the ABNK-2 method.

Moreover, it is calculated from Table \ref{tab:resultEx2} that the speed-up of the ABNK-2 method against the NRK method is around 20, the speed-up of the ABNK-2 method against the MRNK method is around 10, the speed-up of the ABNK-2 method against the MRBNK method is around 300 and the speed-up of the ABNK-2 method against the ABNK-1 method is around 6.

\begin{table}[h]
\centering 
\caption{ Numerical results for Example \ref{exam2}. }\label{tab:resultEx2}  
\resizebox{\textwidth}{!}{  
\begin{tabular}{lcccccccccc} 
\hline
\multirow{2}{*} {$m$}  & \multicolumn{2}{c}{NRK} & \multicolumn{2}{c}{MRNK} & \multicolumn{2}{c}{MRBNK} & \multicolumn{2}{c}{ABNK-1} & \multicolumn{2}{c}{ABNK-2} \\ 
 \cmidrule[0.25mm](lr){6-7} \cmidrule[0.25mm](lr){8-9}\cmidrule[0.25mm](lr){10-11} 
 \cmidrule[0.25mm](lr){2-3} \cmidrule[0.25mm](lr){4-5} \cmidrule[0.25mm](lr){6-7}\cmidrule[0.25mm](lr){8-9} \cmidrule[0.25mm](lr){10-11}  & IT &CPU &IT &CPU &IT &CPU &IT &CPU &IT &CPU  \\ 
\hline
100 &	 219293&	 4.147&	 211476&	 1.425&	 152296&	 68.190&	75059&	 0.718 &	  10464&	 0.181  \\  
200 &	 230374&	 4.837&	 221599&	 1.690&	 153323&	 70.714&	76751&	 1.393 &	  12224&	 0.275  \\  
300 &	 240614&	 5.566&	 231572&	 1.967&	 153979&	 72.977&	78052&	 1.628 &	  11757&	 0.280  \\  
400 &	 251621&	 6.485&	 241693&	 2.601&	 153783&	 92.036&	79352&	 1.818 &	  12312&	 0.341  \\  
500 &	 263013&	 7.184&	 252229&	 2.754&	 154338&	 94.077&	80652&	 1.983 &	  6547&	 0.275  \\  
600 &	 274267&	 8.032&	 262529&	 3.048&	 154690&	 96.766&	82051&	 2.158 &	  12629&	 0.482  \\  
700 &	 284956&	 9.447&	 272705&	 3.548&	 155881&	 100.177&	83357&	 2.316 &	  16631&	 0.678  \\  
800 &	 296740&	 10.020&	 282951&	 3.968&	 155492&	 102.644&	85029&	 2.622 &	  13054&	 0.662  \\  
900 &	 307362&	 14.424&	 293228&	 4.397&	 156521&	 105.705&	86332&	 2.940 &	  13010&	 0.794  \\  
1000 &	 319038&	 13.096&	 303724&	 5.252&	 157102&	 108.886&	87633&	 3.087 &	  13134&	 0.924  \\  
\hline
\end{tabular}  
}   
\begin{tablenotes} 
 \footnotesize               
\item The parameters in the MRBNK, ABNK-1 and ABNK-2 methods are $\theta_1^*=0.5, (\alpha^*,\theta_2^*)=(1.8,0.9)$ and $(\delta^*,\theta_3^*)=(1.0,0.2)$, respectively. 
\end{tablenotes}
 \end{table} 
 
\begin{figure} [!htbp] 
	\centering
	\subfigure[$m=100$]
	{
		\begin{minipage}[t]{0.48\linewidth}
			\centering
			\includegraphics[width=1\textwidth]{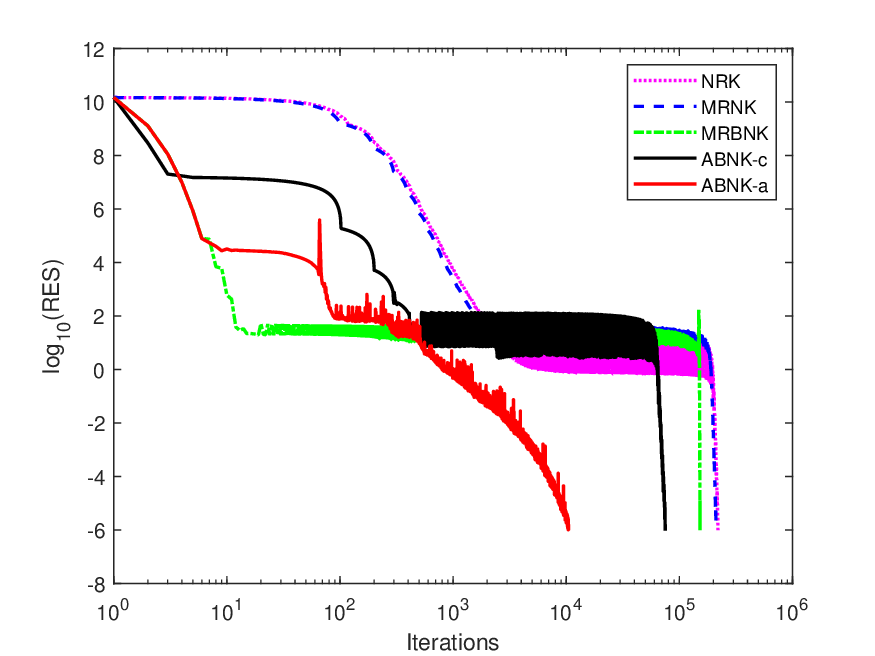}
		\end{minipage}
	}
	\subfigure[$m=100$]
	{
		\begin{minipage}[t]{0.48\linewidth}
			\centering
			\includegraphics[width=1\textwidth]{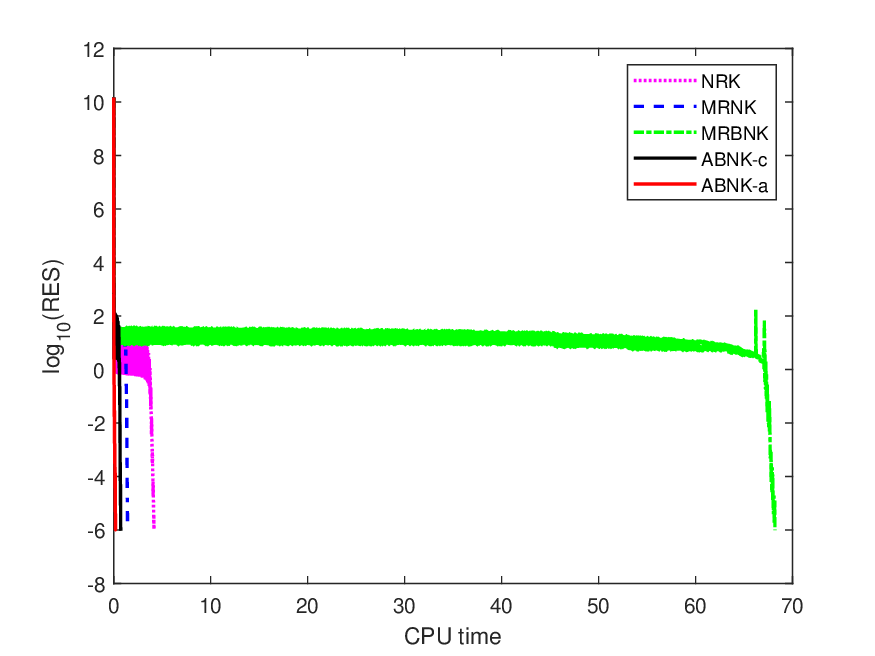}
		\end{minipage}
	}
	\caption{ Curves of residual versus the iteration steps and CPU time in Example \ref{exam2}.} 	
		\label{fig:resvsitcpuEx2}
\end{figure}

In Figure \ref{fig:resvsitcpuEx2}, the curves of the residuals versus the number of iteration steps and the elapsed CPU time for all tested methods are depicted when $m=100$, respectively. 

From Figure \ref{fig:resvsitcpuEx2}, it is observed that the curves of residual versus both the number of iteration steps and the elapsed CPU time for ABNK methods decrease much faster than those of other methods. Moreover, the curves of residual for ABNK-2 method decrease the fastest among all methods. 

\section{Conclusions}\label{secconclu}
 In this paper, a class of averaging block nonlinear Kaczmarz methods is presented for the solution of nonlinear systems of equations. The proposed methods with both constant stepsize and adaptive stepsize converge to the solution of the nonlinear systems. Numerical experiments show the efficiency and robustness of the proposed methods with different stepsizes. Moreover, the proposed method with adaptive stepsize is superior to that with constant stepsize, which outperforms the existing nonlinear Kaczmarz methods both in the aspects of the number of iteration steps and the CPU time.

\bibliographystyle{plain}	 
 
\bibliography{refsabnk.bib}

\begin{thebibliography}{10}

\bibitem{2010BDDKB}
DKR Babajee, MZ~Dauhoo, Mohammad~Taghi Darvishi, A~Karami, and Ali Barati.
\newblock {Analysis of two Chebyshev-like third order methods free from second
  derivatives for solving systems of nonlinear equations}.
\newblock {\em Journal of Computational and Applied Mathematics},
  233(8):2002--2012, 2010.

\bibitem{2015CLS}
Emmanuel~J Candes, Xiaodong Li, and Mahdi Soltanolkotabi.
\newblock {Phase retrieval via Wirtinger flow: Theory and algorithms}.
\newblock {\em IEEE Transactions on Information Theory}, 61(4):1985--2007,
  2015.

\bibitem{2007DB}
Mohammad~Taghi Darvishi and Ali Barati.
\newblock {A third-order Newton-type method to solve systems of nonlinear
  equations}.
\newblock {\em Applied Mathematics and Computation}, 187(2):630--635, 2007.

\bibitem{2020DSS}
Kui Du, Wu-Tao Si, and Xiao-Hui Sun.
\newblock {Randomized extended average block Kaczmarz for solving least
  squares}.
\newblock {\em SIAM Journal on Scientific Computing}, 42(6):A3541--A3559, 2020.

\bibitem{2020FS}
Albert Fannjiang and Thomas Strohmer.
\newblock The numerics of phase retrieval.
\newblock {\em Acta Numerica}, 29:125--228, 2020.

\bibitem{2007HLS}
Markus Haltmeier, Antonio Leitao, and Otmar Scherzer.
\newblock {Kaczmarz methods for regularizing nonlinear ill-posed equations I:
  Convergence analysis}.
\newblock {\em Inverse Problems and Imaging}, 1(2):289, 2007.

\bibitem{2018H}
Wenrui Hao.
\newblock {A homotopy method for parameter estimation of nonlinear differential
  equations with multiple optima}.
\newblock {\em Journal of Scientific Computing}, 74:1314--1324, 2018.

\bibitem{2016K}
Kenji Kawaguchi.
\newblock Deep learning without poor local minima.
\newblock {\em Advances in neural information processing systems}, 29, 2016.

\bibitem{1989L}
Guangye Li.
\newblock {Successive column correction algorithms for solving sparse nonlinear
  systems of equations}.
\newblock {\em Mathematical Programming}, 43(1-3):187--207, 1989.

\bibitem{2022LZB}
Chaoyue Liu, Libin Zhu, and Mikhail Belkin.
\newblock {Loss landscapes and optimization in over-parameterized non-linear
  systems and neural networks}.
\newblock {\em Applied and Computational Harmonic Analysis}, 59:85--116, 2022.

\bibitem{2004MNT}
Kaj Madsen, Hans~Bruun Nielsen, and Ole Tingleff.
\newblock Methods for non-linear least squares problems.
\newblock 2004.

\bibitem{1975M}
Steve~F McCormick.
\newblock {An iterative procedure for the solution of constrained nonlinear
  equations with application to optimization problems}.
\newblock {\em Numerische Mathematik}, 23(5):371--385, 1975.

\bibitem{2022MW}
Cun-Qiang Miao and Wen-Ting Wu.
\newblock {On greedy randomized average block Kaczmarz method for solving large
  linear systems}.
\newblock {\em Journal of Computational and Applied Mathematics}, 413:114372,
  2022.

\bibitem{2019N}
Ion Necoara.
\newblock {Faster randomized block Kaczmarz algorithms}.
\newblock {\em SIAM Journal on Matrix Analysis and Applications},
  40(4):1425--1452, 2019.

\bibitem{2019NRP}
Ion Necoara, Peter Richt{\'a}rik, and Andrei Patrascu.
\newblock {Randomized projection methods for convex feasibility: Conditioning
  and convergence rates}.
\newblock {\em SIAM Journal on Optimization}, 29(4):2814--2852, 2019.

\bibitem{1982PS}
Christopher~C Paige and Michael~A Saunders.
\newblock {LSQR: An algorithm for sparse linear equations and sparse least
  squares}.
\newblock {\em ACM Transactions on Mathematical Software (TOMS)}, 8(1):43--71,
  1982.

\bibitem{2022P}
Kim-Hung Pho.
\newblock {Improvements of the Newton--Raphson method}.
\newblock {\em Journal of Computational and Applied Mathematics}, 408:114106,
  2022.

\bibitem{2022WLB}
Qifeng Wang, Weiguo Li, Wendi Bao, and Xingqi Gao.
\newblock {Nonlinear Kaczmarz algorithms and their convergence}.
\newblock {\em Journal of Computational and Applied Mathematics}, 399:113720,
  2022.

\bibitem{2023XYZ}
A-Qin Xiao, Jun-Feng Yin, and Ning Zheng.
\newblock {On fast greedy block Kaczmarz methods for solving large consistent
  linear systems}.
\newblock {\em Computational and Applied Mathematics}, 42(3):119, 2023.

\bibitem{1995Y}
Tjalling~J Ypma.
\newblock {Historical development of the Newton--Raphson method}.
\newblock {\em SIAM review}, 37(4):531--551, 1995.

\bibitem{2022YLG}
Rui Yuan, Alessandro Lazaric, and Robert~M Gower.
\newblock Sketched newton--raphson.
\newblock {\em SIAM Journal on Optimization}, 32(3):1555--1583, 2022.

\bibitem{2011Y}
Ya-Xiang Yuan.
\newblock {Recent advances in numerical methods for nonlinear equations
  andnonlinear least squares}.
\newblock {\em Numerical algebra, control and optimization}, 1(1):15--34, 2011.

\bibitem{2023ZBLW}
Feiyu Zhang, Wendi Bao, Weiguo Li, and Qin Wang.
\newblock {On sampling Kaczmarz--Motzkin methods for solving large-scale
  nonlinear systems}.
\newblock {\em Computational and Applied Mathematics}, 42(3):126.

\bibitem{2023ZWZ}
Jianhua Zhang, Yuqing Wang, and Jing Zhao.
\newblock {On maximum residual nonlinear Kaczmarz-type algorithms for large
  nonlinear systems of equations}.
\newblock {\em Journal of Computational and Applied Mathematics}, page 115065,
  2023.

\bibitem{2022ZLT}
Yanjun Zhang, Hanyu Li, and Ling Tang.
\newblock {Greedy randomized sampling nonlinear Kaczmarz methods}.
\newblock {\em arXiv preprint arXiv:2209.06082}, 2022.

\end{thebibliography}

\end{document}